\documentclass[leqno]{article}
\usepackage{latexsym,amssymb,amsmath}
\title{\bf Tame coverings of arithmetic schemes}
\author{by Alexander Schmidt\\ at Heidelberg}
\date{}
\newtheorem{itheorem}{Theorem}
\newtheorem{theorem}{Theorem}[section]
\newtheorem{example}[theorem]{Example}

\newtheorem{remark}[theorem]{Remark}
\newtheorem{prop}[theorem]{Proposition}
\newtheorem{corol}[theorem]{Corollary}
\newtheorem{defi}[theorem]{Definition}
\newtheorem{lemma}[theorem]{Lemma}

\def\demo #1{
    \vskip-\lastskip \vskip 12pt plus 3pt minus 9pt%
\noindent{\bf #1}\enspace}
\def\noproof{{\unskip\nobreak\hfill\penalty50\hskip2em\hbox{}%
     \nobreak\hfill$\square$\parfillskip=0pt%
     \finalhyphendemerits=0\par}}
\def\enddemo{\ifmmode\eqno\square\else\noproof\vskip 12pt plus 3pt minus 9pt
\fi}
\def\diagram{\renewcommand\arraystretch{1.5} $$ \begin{array}}
\def\enddiagram{\end{array} $$ \renewcommand\arraystretch{1}}

\def\rmnewname#1{\expandafter\gdef\csname#1\endcsname{{\mathop{\rm
#1}\nolimits}}}
\def\itnewname#1{{\expandafter\gdef\csname#1\endcsname{{\mathop{\it
#1}\nolimits}}}}

 \itnewname{ab}
 \itnewname{can}
 \rmnewname{codim}
 \rmnewname{cone}
 \itnewname{cor}
 \rmnewname{CH}
 \rmnewname{div}
 \rmnewname{Div}
 \itnewname{et}
 \rmnewname{Ext}
 \rmnewname{Frob}
 \rmnewname{Hom}
 \itnewname{id}
 \rmnewname{Ker}
 \rmnewname{length}
 \rmnewname{Mor}
 \rmnewname{Pic}
 \rmnewname{Quot}
 \itnewname{rec}
 \itnewname{res}
 \rmnewname{Sch}
 \itnewname{sh}
 \rmnewname{Sm}
 \rmnewname{SmCor}
 \itnewname{sp}
 \rmnewname{Spec}
 \rmnewname{supp}
 \itnewname{tame}
 \rmnewname{top}
 \rmnewname{Tor}
 \itnewname{Zar}

 \def\gP{{\mathfrak P}}
 \def\p{{\mathfrak p}}

 \def\Q{{\mathbb Q}}

 \def\Z{{\mathbb Z}}

 \def\noi{\noindent}
 \def\bs{\bigskip}
 \def\ms{\medskip}

 \def\Cal#1{{\cal #1}}

 \def\lang{\longrightarrow}
 \def\longeq{\hbox{$=\!=$}}
 \def\mapr#1{\mathrel{\mbox{$\stackrel{#1}{\longrightarrow}$}}}
 
 \def\mapd#1{\Big\downarrow\rlap{$\vcenter{\hbox{$\scriptstyle{{#1}}$}}$}}
 
 \def\liso{\mathrel{\hbox{$\longrightarrow$} \kern-14pt\lower-4pt%
 \hbox{$\scriptstyle\sim$}\kern7pt}}
 \def\r@iso{\mathrel{\lower2pt\hbox{$\scriptstyle\sim$}%
 \kern-8pt\hbox{$\rightarrow$}}}
 \def\riso#1{\mathrel{\stackrel{\!\!#1}{\r@iso}}}
 \def\lr@iso{\mathrel{\kern6pt\lower2pt\hbox{$\scriptstyle\sim$}%
  \kern-12pt\hbox{$\longrightarrow$}}}
 \def\lriso#1{\mathrel{\mathop{\lr@iso}\limits^{#1}}}
 \def\eqd{\Big\|}

 \def\mapd#1{\Big\downarrow\rlap {$\vcenter{\hbox{$\scriptstyle{{#1}}$}}$}}
 \def\surjr#1{\stackrel{#1}{\hbox{$\relbar \! \! \! \twoheadrightarrow$}}}
 
 \def\injr#1{\stackrel{#1}{\lhook\joinrel\relbar\joinrel\rightarrow}}
 
 \def\surjd#1{\lower4pt\hbox{$\downarrow$}\kern-5.75pt\Big\downarrow\rlap%
  {$\vcenter{\hbox{$\scriptstyle{{#1}}$}}$}}
 \def\surju#1{\lower-4pt\hbox{$\uparrow$}\kern-5.75pt\Big\uparrow\rlap%
  {$\vcenter{\hbox{$\scriptstyle{{#1}}$}}$}}

 \def\inju#1{
 \setlength{\unitlength}{0.1pt}
 \begin{picture}(40,120)
 \put(20,-21){\vector(0,1){142}} \put(-16.3,-50){\mbox{$\scriptscriptstyle
 \cup $}} \put(47,10){\mbox{$\scriptstyle #1 $}}
 \end{picture}}

 \def\projlim{\mathop{\lim\limits_{{\longleftarrow}}\,}}

\def\tcup{\mathop{\mathord{\cup}\mkern-8.5mu^{\hbox{.}}\mkern 9mu}}
\def\ttcup{\mathop{\mathord{\bigcup}\mkern-10mu^{\hbox{.}}\mkern 9mu}}

\def\TTcup{\textstyle  \ttcup\limits}
\def\Tcup{\textstyle \tcup}

 \def\hang{\hangindent\Itemindent}
 \def\textindent#1{\hskip\Itemindent\llap{\hbox{\rm #1}\enspace}\ignorespaces}
 \def\Item{\par\noindent\hang\textindent}
 \def\ItemItem{\par\noindent\hskip\Itemindent \hangindent2\Itemindent \textindent}
 \newdimen\Itemindent \Itemindent=.9cm

\begin{document}
\maketitle

The objective of this paper is to investigate tame fundamental groups of schemes of finite
type over $\Spec (\Z)$. More precisely, let $X$ be a connected scheme of finite type over
$\Spec(\Z)$ and let $\bar X$ be a compactification of $X$, i.e.\ a scheme which is proper and
of finite type over $\Spec(\Z)$ and which contains $X$ as a dense open subscheme. Then the
tame fundamental group of $X$ classifies finite \'{e}tale coverings of $X$ which are tamely
ramified along the boundary $\bar X -X$, in particular, the tame fundamental group
$\pi_1^t(\bar X,\bar X -X)$ is a quotient of the \'{e}tale fundamental group $\pi_1(X)$. Our
interest in the tame fundamental group arises from the observation that it seems to be the
maximal quotient of $\pi_1(X)$ which is \lq visible\rq\ via class field theory by algebraic
cycle theories (see \cite{S-S1}, \cite{S} and \cite{S-S2} for more precise statements on \lq
tame class field theory\rq).

\ms Coverings of a regular scheme which are tamely ramified along a normal crossing divisor
have been studied in \cite{SGA1}, \cite{G-M}. In this paper we consider tame ramification
along an arbitrary Zariski-closed subset. The main reason for this is the lack of good
desingularization theorems in positive and mixed characteristics. A simple imitation of the
definition of tame ramification in the normal crossing case proves to be not useful in the
more general situation. We give a definition of tameness in the general situation in
section~\ref{tamesect} and we show that it coincides with the previous one in the normal
crossing case. For Galois coverings of normal schemes our definition of tameness coincides
with that proposed by Abbes  \cite{Ab} and with the notion of ``numerical tameness'' defined
by Chinburg/Erez \cite{C-E}.

 Naturally, the question arises whether the tame
fundamental group $\pi_1^t(\bar X, \bar X-X)$ is independent of the choice of the
compactification $\bar X$. At the moment, we can answer this question only for the maximal
pro-nilpotent quotient. In section~\ref{valsect}, we consider discrete valuations of higher
rank associated to Parshin chains. We investigate the connection between the tameness of a
covering and the behaviour of the associated higher dimensional henselian fields and we show
that a finite covering of a regular arithmetic scheme with nilpotent Galois group is tamely
ramified if and only if all associated extensions of higher dimensional henselian fields are
tame. Then we use this fact in the proof of the following

\begin{itheorem}\label{i1}
Let $X$ be a regular and connected scheme of finite type over $\Spec(\Z)$ and assume that
there exists a regular scheme $\bar X$ which is proper over $\Spec(\Z)$ and which contains
$X$ as an open subscheme. Then $\pi_1^t(\bar X,\bar X-X)^{\hbox{\it\scriptsize pro-nil}}$
depends only on $X$ and not on the choice of\/ $\bar X$.
\end{itheorem}

 If $X$ is a smooth quasiprojective variety over a finite field $\mathbb
F$, then its abelianized tame fundamental group $\pi_1^t(X)^{ab}=\pi_1^t(\bar X,\bar X
-X)^{ab}$ is infinite, but its degree-zero part, i.e.\ the kernel of the natural homomorphism
$$
\pi_1^t(X)^{ab} \longrightarrow G(\bar {\mathbb F}\,|\,{\mathbb F}) \cong \hat \Z
$$

\noindent
is finite (see \cite{S-S1}). In the arithmetic case, i.e.\ when $X$ is flat over $\Spec(\Z)$,
we show the following theorem in which $\bar X$ is not assumed to be proper.

\begin{itheorem} \label{introfin}
Let ${\Cal O}$ be the ring of integers in a finite extension $k$ of \/ $\mathbb Q$. Let $\bar
X$ be a flat $\Cal O$-scheme of finite type whose geometric generic fibre $\bar X
\otimes_{\Cal O}\bar k$ is connected. Assume that $\bar X$ is normal and that the morphism
$\bar X \to \Spec (\Cal O)$ is surjective. Let $X$ be an open subscheme of $\bar X$. Then the
abelianized tame fundamental group $\pi_1^t(\bar X,\bar X-X)^{ab}$ is finite.
\end{itheorem}
In the special case $X=\bar X$, we obtain the finiteness of the \'{e}tale fundamental group
$\pi_1(\bar X)^{ab}$ which sharpens a theorem of Katz and Lang (\cite{K-L}, th.3) by
weakening the assumption \lq smooth\rq\ to \lq normal\rq.

\section{Tame coverings} \label{tamesect}

 The concept of tame ramification stems from number theory: A
finite extension of number fields $L|K$ is called tamely ramified at a prime $\gP$ of $L$ if
the associated extension of completions $L_\gP|K_\gP$ is a tamely ramified extension of local
fields. The latter means that the ramification index is prime to the characteristic of the
residue field. It is a classical result that composites and towers of tamely ramified
extensions are again tamely ramified. This concept generalizes to arbitrary discrete
valuation rings by requiring that the associated residue field extensions are separable.

\ms As is well known, the concept of unramified extensions has found its generalization to
arbitrary schemes by the notion of \'{e}tale coverings. Tame ramification along a normal
crossing divisor of a regular scheme has been studied in \cite{SGA1}, \cite{G-M}. Let us
first recall the definition given there. We assume throughout this paper that all schemes
under consideration are noetherian.

\ms Let $X$ be a regular and connected scheme and let $D=\sum
D_i$ be a divisor on $X$. We say that $D$ has {\bf normal crossings}, if \'{e}tale locally around
every point $x \in \supp (D) \subset X$, we have
\[
D_i= \sum_j \div (s_{ij}),
\]
where $(s_{ij})_{i,j} \in \Cal O_{X,x}$ is part of a regular system of parameters of the
regular local ring $\Cal O_{X,x}$.

\ms Let $X$ be regular and connected, $D\subset X$ a divisor with
normal crossings and  $U=X-D$. Let $P_i$, $i=1,\ldots,n$, be the generic points of the
irreducible components $D_1, \ldots, D_n$ of $D$. Then the local rings $\Cal O_{X,P_i}$ are
discrete valuation rings, inducing discrete valuations $v_1,\ldots,v_n$ on the function field
of $X$. Let $U'\to U$ be a finite \'{e}tale morphism and assume for simplicity that $U'$ is
connected. Let $X'$ be the normalization of $X$ in the function field $K'$ of $U'$.
\begin{defi}{\rm (\cite{G-M}, 2.2.2.)} \label{tame1}
The finite \'{e}tale covering $U' \to U$ is called {\bf tamely ramified} (or tame, for short)
along $D$ if the extension of function fields $K'|K$ is tamely ramified at the discrete
valuations associated to $D_1, \ldots, D_n$.
\end{defi}
\'{E}tale locally, tame coverings are of a very simple structure by the following theorem which
is known under the name Generalized Abhyankar's Lemma.

\begin{theorem} {\rm (\cite{SGA1}, Exp.\ XIII, 5.3.0)}
Let $X$ be a strictly henselian local regular scheme of residue characteristic $p>0$, $D=
\sum_{i=1}^r \div (f_i)$ a divisor with normal crossings on $X$ and $U= X-D$. Then every
connected finite \'{e}tale covering of $U$ which is tamely ramified along $D$ is a quotient of a
(tamely ramified) covering of the form \label{abhy}
\[
U'= U[T_1, \ldots, T_r]/(T_1^{n_1}-f_1, \ldots, T_r^{n_r}-f_r),
\]
where the $n_i$ are natural numbers prime to $p$.
\end{theorem}

Tame coverings as defined above satisfy the axioms of a Galois theory, so (omitting the base
point) there exists a profinite group $\pi_1^t(X,D)$, which is a quotient of $\pi_1(U)$ and
which classifies \'{e}tale coverings of $U$ which are tame along $D$.

\ms If $X$ is strictly local and $D= \sum_{i=1}^r \div (f_i)$ a
normal crossing divisor, then theorem \ref{abhy} yields a natural isomorphism
\[
\pi_1^t(X,D) \cong \Big(\prod_{\ell \not = p} \Z_\ell(1)\Big)^r \quad \stackrel{notation}{=}
\left(\hat \Z^{(p')}(1)\right)^r .
\]

\bs If, more generally, $D \subset X$ is an arbitrary divisor on
$X$, then (as already remarked in \cite{G-M}, 2.2.3.4, without further elaboration) the above
definition of tame ramification is not the `correct' one. For example, it is not stable under
base change.

\begin{example}\label{counterex}\rm
Let $X=\Spec(\Z[T])$ be the affine line over $\Spec(\Z)$ and consider the divisor
\[
D=\div (T+4) + \div(T-4),
\]
which is not a normal crossing divisor. Let $K=\Q(T)$ be the function field of $X$ and
$U=X-D$. Put $f=(T+4)(T-4)=T^2-16$, $L=K(\sqrt{f})$ and consider the normalization $X_L$ of
$X$ in $L$. The ramification locus of $X_L \to X$ is either $D$ or $D \cup X_2$, where $X_2$
is the unique vertical divisor on $X$ over characteristic~$2$. Let us show that $X_L \to X$
is unramified at the generic point of $X_2$. This is equivalent to the statement that $L|K$
is unramified at the unique discrete valuation $v_2$ of $K$ which corresponds to the prime
ideal $2\Z[T]\subset \Z[T]$. Therefore it suffices to show that $f$ is a square in the
completion $K_2$ of $K$ with respect to $v_2$. Consider the polynomial $F(X)=X^2
-f=X^2-T^2+16$. We have  $F(T)\equiv 0 \bmod 16$ and the derivative $F'(T)= 2 T$ has the
exact $2$-valuation $1$. By the usual approximation process (cf.\ \cite{Se} 2.2.\ th.1), we
see that $f$ has a square root in $K_2$.  Hence the ramification locus of $X_L \to X $ is
exactly $D$, and since $D$ is the sum of horizontal prime divisors, the morphism $U_L \to U$
is tame along $D$ (in the naive sense).

Now consider the closed subscheme $Y \subset X$ given by the equation \hbox{$T=0$}, so $Y
\cong \Spec(\Z)$. Then $D_Y= D \cap Y$ is the point  on $Y$ which corresponds to the prime
number $2$. Let $V= U \cap Y= Y -D_Y$. The base change $V'=U_L \times_U V \to V$ is the
normalization of $V\cong \Spec(\Z[\frac{1}{2}])$ in  $\Q(\sqrt{-1})$.  But $2$ is wildly
ramified in $\Q(\sqrt{-1})$, and so $V' \to V$ is not tame along  $D_Y$.
\end{example}

\bigskip
Let us now give a definition of tame ramification in the general situation which will be
shown in proposition~\ref{coincide} to generalize definition~\ref{tame1}. Let $X$ be a
scheme, $Y \subset X$ a closed subscheme and $U=X-Y$ the open complement. For a point $y \in
Y$ we write $X_y^{\sh}$ for $\Spec(\Cal O_{X,y}^{\sh})$, where $\sh$ means strict
henselization. By abuse of notation, we write $U_y^{\sh}$ for the base change $U \times_X
X_y^{\sh}$. The scheme $U_y^{\sh}$ is empty if $y \notin \bar U$.

\begin{defi} \label{tame2}
We say that a finite \'{e}tale covering $U' \to U$ is {\bf tamely ramified} along $Y$ if for
every point $y \in Y$ such that $U_y^{\sh}$ is nonempty, the base change
\[
U' \times_U U_y^{\sh} \lang U_y^{\sh}
\]
can be dominated by an \'{e}tale covering of the form
\[
V_1 \Tcup \cdots \Tcup V_r \lang U_{y}^{\sh},
\]
such that each $V_i$ is a connected \'{e}tale Galois covering of its image in $U_y^{\sh}$ and the
degree of\/ $V_i$ over its image is prime to the characteristic of\/ $k(y)$.
\end{defi}

The next lemma follows in a straightforward manner from the definition of tame ramification.

\begin{lemma}\label{tamebasechange} Let $X$ be a scheme, $U \subset X$ an
open subscheme and $Y=X-U$. Let $f: X_1 \to X$ be a morphism of schemes, $U_1=f^{-1}(U)$ and
$Y_1=f^{-1}(Y)$. If\/ $U' \to U$ is a finite \'{e}tale morphism which is tamely ramified along\/
$Y$, then the base change
\[
U_1'= U' \times_U U_1 \lang U_1
\]
is tamely ramified along\/ $Y_1$.
\end{lemma}

The notion of a finite \'{e}tale morphism $U' \to U$ which is tamely ramified along $Y=X-U$ is
independent from $X$ in the following sense.

\begin{lemma} \label{indep}
Suppose that $U$ is contained as an open subscheme in schemes $\widetilde X$ and $X$ with
closed complements $\widetilde Y$ and $Y$, respectively. Assume that there exists a finite
morphism $\pi: \widetilde X \to X$ making the diagram \diagram{ccc}
 U&\injr{}&\widetilde X\\
 \eqd && \mapd{\pi} \\
 U& \injr{}&X
\enddiagram
commutative. Then a finite \'{e}tale morphism $U' \to U$ is tamely ramified along $\widetilde Y$
if and only if it is tamely ramified along $Y$.
\end{lemma}

\demo{Proof:} Let $y \in Y$ be a point and consider the cartesian diagram \diagram{ccc}
 \widetilde X \times_{X} X_{y}^{\sh}&\lang & X_{y}^{\sh} \\
 \inju{}&& \inju{}\\
U \times_{X} X_{y}^{\sh}& \longeq &U_{y}
\enddiagram
Since $\pi: \widetilde X \to X$ is finite, $ \widetilde X \times_{X} X_{y}^{\sh}$ is strictly
henselian, and
\[
 \widetilde X \times_{X} X_{y}^{\sh} \cong \TTcup_{i} X_{\tilde y_{i}}^{\sh} \;,
\]
where the $\tilde y_{i}$ are the finitely many points of $\widetilde X$ lying above $y$.
Therefore, $\pi$ induces a natural isomorphism
\[
\TTcup_i U_{\tilde y_{i}}^{\sh} \liso U_{y}^{\sh}.
\]
Now the statement of the lemma follows easily from the definition of tame ramification.
\enddemo

\begin{remark}\rm
 Assume that $X$ is excellent or, slightly weaker, that all local rings
of $X$ are Nagata rings. Then there exists a unique maximal scheme $\widetilde X$ satisfying
the conditions of lemma \ref{indep}: the normalization of $X$ in $U$. This scheme is
constructed in the following way. Denoting the open immersion by $j: U \hookrightarrow X$,
the sheaf $j_*(\Cal O_U)$ is quasicoherent and contains $\Cal O_X$. The integral closure of
$\Cal A$ of $\Cal O_X$ in $j_*(\Cal O_U)$ is then a coherent sheaf of $\Cal O_X$-algebras and
the normalization of $X$ in $U$ is defined as $\Cal Spec(\Cal A)$ (cf.\ \cite{Mi}, I, \S1,
proof of th.\ 1.8).
\end{remark}

If $U' \to U$ is a tame covering, then, by Zariski's main theorem, the morphism $U' \to X$
factors in the form $U' \mapr{j} X' \mapr{\pi} X$ with $\pi$ finite and $j$ an open
embedding. In general, $X'$ is not unique but if $X_1'$ and $X_2'$ are such schemes, we find
a third scheme $X_3'$ dominating $X_1'$ and $X_2'$. If $X$ is normal and connected, then the
normalization of $X$ in the function field of $U'$ is maximal among the schemes which are
finite over $X$ and contain $U'$ as an open dense subscheme. The following lemma follows
easily by the same considerations as in the proof of lemma \ref{indep}.

\begin{lemma}
Let
$$
X_1\, \mapr{\pi_1} \, X_2\,  \mapr{\pi_2}\, X_3
$$
be finite surjective morphisms, let $U_3 \subset X_3$ be an open subscheme and denote by
$U_2$ and $U_1$ the corresponding inverse images. Suppose that $\pi_1|_{U_1}: U_1\to U_2$ and
$\pi_2|_{U_2}: U_2 \to U_3$ are \'{e}tale. Furthermore, let\/ $Y_i$, $1 \leq i \leq 3$, be the
closed complement of\/ $U_i$ in $X_i$.

\ms Then $(\pi_2 \circ \pi_1)|_{U_1}: U_1 \to U_3$ is tamely
ramified along $Y_3$ if and only if $\pi_1|_{U_1}: U_1 \to U_2$ is tamely ramified along
$Y_2$ and $\pi_2|_{U_2}: U_2 \to U_3$ is tamely ramified along $Y_3$.
\end{lemma}

Consider the category ${\bf FEtT}/(X,Y)$ of $U$-schemes, finite and \'{e}tale over $U$ with tame
ramification along $Y$. It is a full subcategory of the category  ${\bf FEt}/U$ of
$U$-schemes, finite and \'{e}tale over $U$. If\/ $\bar{x} \rightarrow U$ is a geometric point of
$U$, we consider the fibre functor
\[
{\bf F}:\quad {\bf FEtT}/(X,Y) \mapr{} (Sets); \quad (V\to U) \longmapsto \Mor_U(\bar x,V),
\]
which is just the restriction of the usual fibre functor induced by $\bar x$ on ${\bf FEt}/U$
to ${\bf FEtT}/(X,Y)$. One observes that the category ${\bf FEtT}/(X,Y)$ together with this
fibre functor satisfies the axiomatic conditions for a Galois theory (see \cite{SGA1}
Exp.V,4). Thus one obtains a profinite fundamental group
\[
\pi_1^t(X,Y;\bar x)
\]
which we call the {\bf tame fundamental group} of $U$ with base point $\bar x$ and with
respect to the embedding $U \hookrightarrow X$. It is a quotient of the \'{e}tale fundamental
group $\pi_1(U;\bar x)$. If $\bar x' \hookrightarrow U$ is another base point in the same
connected component of $U$, then the fundamental group $\pi_1^t(X,Y;\bar x')$ is isomorphic
to $\pi_1^t(X,Y;\bar x)$, the isomorphism is determined up to an inner automorphism. In the
following, $U$ will always be connected and we will omit the base point from the notation.
If $X$ is a scheme of characteristic $0$ (i.e.\ all residue fields have characteristic~$0$),
then
\[
\pi_1^t(X,Y)=\pi_1(U).
\]
If $X$ is normal and connected, then the functor ${\bf FEt}/X \to {\bf FEt}/U$ is fully
faithful, and we obtain surjections
\[
\pi_1(U) \twoheadrightarrow \pi_1^t(X,Y) \twoheadrightarrow \pi_1(X).
\]
The next proposition follows in a straightforward manner from the theorem on the purity of
the branch locus of Zariski-Nagata (\cite{SGA2}, Exp.X, th.\ 3.4).

\begin{prop}
Assume that $X$ is regular and connected, and that $Y$ is of\/ co\-di\-men\-sion $\geq 2$ in
$X$. Then we have natural isomorphisms
\[
\pi_1(U) \liso \pi_1^t(X,Y) \liso \pi_1(X).
\]
\end{prop}

Let $U' \to U$ be a finite \'{e}tale morphism. The {\bf wild locus} $W_{(U'\to U)}$, i.e.\ the
set of points $y \in Y$ such that $U' \to U$ is not tamely ramified at $y$, is a closed
subscheme of\/ $Y$. The theorem of Zariski-Nagata on the purity of the branch locus says
that if $X$ is regular, the ramification locus of a quasifinite ´dominant morphism to $X$ is
always pure of codimension $1$ in $X$. The same is not true of the wild locus, since in
example~\ref{counterex} we constructed a cyclic cover of a regular scheme of dimension $2$
whose wild locus consists of a single closed point, i.e.\ is of codimension $2$. But if $X$
is equicharacteristic, we have the

\begin{prop} \label{codim1}
Assume that $X$ is regular and connected and that all points of $X$ have the same residue
characteristic. Let $U' \to U$ be a finite nilpotent \'{e}tale covering.  Then $W_{(U'\to U)}$ is
either empty or pure of codimension~$1$ in $X$.
\end{prop}

\demo{Proof:} We may assume that all points on $X$ have the common residue characteristic
$p>0$, because otherwise the wild locus is empty. Let $K$ be the function field of $X$ and
$L$ the function field of $U'$. Then $L|K$ is a finite nilpotent extension and $U'=U_L$ is
the normalization of $U$ in $L$.  Writing $L$ as a composite $L=L_1L_2$, where $[L_1:K]$ is
prime to $p$ and $[L_2:K]$ is a $p$-power, then  $U_{L_1} \to U$ {\it is} tamely ramified
along $Y$. The ramification locus of $X_{L_2} \to X$ is the same as the wild locus
$W_{U_{L_2}\to U}$ since $G(L_2|K)$ is a $p$-group and all points of $X$ have residue
characteristic $p$. Since $X$ is regular, the ramification locus of $X_{L_2} \to X$ is either
empty or pure of codimension $1$ and coincides with $W_{(U_{L_2}\to U)}=W_{(U'\to U)}$.
\enddemo

\begin{remark}\rm
Under the assumptions of proposition~\ref{codim1}, let $D_1, \ldots, D_r$ be the irreducible
components of $Y=X-U$ which are of codimension~$1$ in $X$. Then a finite nilpotent covering
$U' \to U$ is tamely ramified along $Y$ if and only if the discrete valuations $v_1,
\ldots,v_r$ of $K(X)=K(U)$ which are associated to $D_1, \ldots, D_r$ are tamely ramified in
$K(U')|K(U)$. This justifies the \lq naive\rq\ definition of tame ramification used in
\cite{S-S1}.
\end{remark}

>From now on assume that $X$ is a normal connected scheme and we denote the function field of
$X$ by $K$. Every connected finite \'{e}tale covering $\widetilde X \to X$ coincides with the
normalization of $X$ in some finite separable extension $L$ of $K$. Therefore (omitting
suitable chosen base points from the notation), we obtain a natural surjection $ G(\bar K |K)
=\pi_1(\Spec(K)) \surjr{} \pi_1(X)$ and an isomorphism
\[
G(\bar K|K) = \projlim_{U \subset X} \pi_1(U),
\]
where $U$ runs through the open subschemes of $X$.

\ms Let us collect some facts on decomposition and inertia of
integrally closed domains (see \cite{B-Comm} Ch.\ V, \S\S2,3). Let $A$ be an integrally
closed domain with quotient field $K$, $L$ a Galois extension of $K$ and $B$ the integral
closure of $A$ in $L$. Let $\gP$ be a prime ideal of $B$ and $\p=\gP \cap A$ the prime ideal
in $A$ lying under $\gP$. Let $k(\gP)=B_\gP/\gP B_\gP$ and $k(\p)=A_\p/\p A_\p$ be the
residue fields of $\gP$ and $\p$, respectively. Let $G=G(L|K)$ be the Galois group. Then $G$
acts transitively on the set of prime ideals of $B$ lying over $\p$. Then one has the
following subgroups in the Galois group.
\begin{description}
\item{-} $Z=Z_\gP(L|K)=\{\sigma \in G\, | \, \sigma(\gP)=\gP \}$ - the decomposition
group,
\item{-} $T=T_\gP(L|K)=\{\sigma \in G_\gP \, |\, \bar{\sigma}=\id: B/\gP \to B/\gP
\}$ - the inertia group.
\end{description}
$T$ is a normal subgroup in $Z$. Let $A^T$ and $A^Z$ be the integral closures of $A$ in $L^T$
and $L^Z$ respectively. Put $\gP^T= \gP \cap A^T$ and $\gP^Z = \gP \cap A^Z$. For a proof of
the following proposition see \cite{Ra}, X, th.2.

\begin{prop} If $L= \bar K$ is the separable closure of
$K$, then \ms
 \Item{\rm (i)} $\gP$ is the only prime ideal in $B$ extending the prime
 ideal $\gP^Z$ in $A^Z$, \ms
 \Item{\rm (ii)} $(A^Z_{\gP^Z},\gP^Z)$ is the
 henselization of $(A_\p,\p)$, in particular, $k(\gP^Z)=k(\p)$, \ms
 \Item{\rm (iii)} $(A^T_{\gP^T},\gP^T)$ is the strict henselization of $(A_\p,\p)$,
in particular, $k(\gP^T)$\\ is the separable closure of $k(\p)$.
\end{prop}

Assume that $X$ is normal and connected.  Then the normalization of $X$ in a finite separable
extension of its function field is finite over $X$. As is well known, for every point $x$ on
$X$, the scheme $X_x^{\sh}$ is noetherian, normal and connected.

\begin{corol} \label{tame-inert}
Let $X$ be a normal connected scheme, $U\subset X$ an open subscheme and\/ $Y=X-U$ the closed
complement. Let $L$ be a finite Galois extension of the function field $K$ of $X$ and let
$U_L$, $X_L$ be the normalizations of $U$ and $X$ in $L$. Then the following are equivalent

\bs { \Item{\rm (i)} $U_L \to U$ is \'{e}tale and tamely ramified
along $Y$, \ms \Item{\rm (ii)} For each point $\gP \in X_L$ the following holds, where $p$
 denotes the residue characteristic of\/ $\gP$: \smallskip
   \ItemItem{\rm a)} {$T_\gP(L|K)$ is of order prime to $p$ }\smallskip
   \ItemItem{\rm b)} {$T_\gP(L|K)=0$ if\/ $\gP \in U_L$.}

  }
\end{corol}

\demo{Proof:} The only nontrivial point is to see that $U_L \to U$ is flat, but this is well
known (cf.\ \cite{Mi}, I,\S3, th.3.21).
\enddemo

\begin{prop}\label{coincide}
Let $X$ be a regular connected scheme and $D \subset X$ a normal crossing divisor. Then both
notions of tame ramification coincide.
\end{prop}

\demo{Proof:} Let (notation as above) $U_L \to U$ be tame in the sense of definition
\ref{tame2} and assume that $v$ is a discrete valuation on $K$ which corresponds to the
generic point $P$ of an irreducible component of $D$. Let (for any embedding to the separable
closure of $K$) $K_P^{\sh}$ be the quotient field of the strict henselization of the discrete
valuation ring $\Cal O_{X,P}$. Then, by definition, $LK_P^{\sh}$ is a product of fields each
of which is contained in a finite Galois extension of degree prime to $p= \hbox{\rm char} \,
k(P)$ of $K_P^{\sh}$. Therefore $v$ is tamely ramified in $L|K$. On the other hand, if all
these discrete valuations are tamely ramified in $L|K$, then it follows from the Generalized
Abhyankar's lemma (theorem \ref{abhy}) that $U_L \to U$ is tamely ramified in the sense of
definition \ref{tame2}.\enddemo

\ms\noi {\bf Remark:} Corollary~\ref{tame-inert} asserts that for coverings of normal schemes
our definition of tameness coincides with that given in \cite{Ab}. Furthermore (in the
situation of corollary~\ref{tame-inert}), if $U_L \to U$ is tame, then $X_L \to X$ is
numerically tame in the sense of \cite{C-E}. It is also not difficult to show  the inverse
implication.

\section{Valuations and tame ramification} \label{valsect}

The aim of this section is to relate tame covers of schemes of finite type over $\Spec(\Z)$
to discrete valuations of higher rank. The quotient fields of the henselizations of such
valuation rings are called higher dimensional henselian fields and are the basic constituents
of Kato's and Saito's higher dimensional class field theory \cite{K-S}.

\ms Let us recall some facts from valuation theory (see
\cite{B-Comm}, \cite{Z-S}, Ch.\ VI). Let $K$ be a field. A subring $V \subset K$ is called
valuation ring if $ x \in K - V\Longrightarrow x^{-1} \in V$. A valuation ring is local with
maximal ideal ${\mathfrak m}_V = \{x \in V\,|\, x^{-1} \not \in V\} $ and is integrally
closed in $K$. The quotient $ \Gamma_V= K^\times / V^\times $ is a totally ordered abelian
group ($\bar{x} \leq \bar {y} \Leftrightarrow x^{-1}y \in V$) and is called the value group
of $V$. Usually, one extends the natural map $ v: K^\times \lang \Gamma_V, $ to $K$ by
setting $v(0)=\infty$. As is well known, the valuation ring $V$ represents an equivalence
class of abstractly defined valuations of $K$; the valuation $v$ is a natural representative
of this class.

If $L|K$ is a finite separable field extension, then there are at least one and at most
finitely many valuations $w$ of $L$ extending $v$, i.e.\ such that $W \cap K =V$, where $W
\subset L$ is the valuation ring of $w$. We have an induced injective homomorphism $\Gamma_V
\hookrightarrow \Gamma_W$. The index of $\Gamma_V$ in $\Gamma_W$ is called the ramification
index and is usually denoted by $e=e_{w|v}$. We also have an associated field extension
$k(v)=V/{\mathfrak m}_V \hookrightarrow k(w)=W/{\mathfrak m}_W$ whose degree is denoted by
$f=f_{w|v}$. The inequality
\begin{equation}\label{defectineq}
\sum_{w|v} e_{w|v} f_{w|v} \leq [L:K],
\end{equation}
(see \cite{B-Comm} Ch.6, \S8.3. th.1) shows, in particular, that all these numbers are
finite. $v$ is called {\bf defectless} in $L$ if equality holds in (\ref{defectineq}). By
\cite{B-Comm}, VI, \S8.5, cor.2, a discrete valuation is defectless in a finite separable
extension.

\ms Let $L$ be a Galois extension of $K$, $W$ a valuation ring of $L$ and
$V=W\cap K$. Let $\gP$ be the unique maximal ideal of $W$. Then, besides the inertia and
decomposition group, we have the {\bf ramification group}
\[
 R_\gP= R_\gP(L|K)=\{ \sigma \in T_\gP(L|K)\,|\, \frac{\sigma x}{x}
\equiv 1 \bmod \gP \hbox{ for all } 0\not = x \in L\}.
\]
A proof of the following proposition can be found in \cite{End}, \S20.
\begin{prop} \label{tameprop}
If the residue characteristic of\/ $V$ is zero, then $R_\gP=0$ and $T_\gP$ is abelian. If the
residue characteristic of\/ $V$ is $p>0$, then

\ms \Item{\rm (i)} $R_\gP$ is a normal subgroup in
$Z_\gP$,\smallskip \Item{\rm (ii)} $R_\gP$ is the unique $p$-Sylow subgroup in
$T_\gP$,\smallskip \Item{\rm (iii)} $T_\gP/R_\gP$ is an abelian group of order prime to $p$.
\end{prop}

The decomposition, inertia, ramification groups of different valuation rings $W$ of $L$
extending $V$ are conjugate in $G(L|K)$. One says that $W|V$ (resp.\ $w|v$) is {\bf tamely
ramified} if the ramification group is trivial. $V$ is tamely ramified in $L|K$ if $W|V$ is
tamely ramified for one (every) valuation ring $W$ of $L$ extending $V$. One says that $V$ is
tamely ramified in a separable extension if it is tamely ramified in the Galois closure of
this extension.

If $V$ is tamely ramified in $L|K$, then it is defectless in $L|K$ (see \cite{End}, cor.
20.22). A valuation $v$ on $K$ is called {\bf\boldmath discrete valuation of rank $n$} if its
value group is (as an ordered group) isomorphic to $\Z\times \cdots \times \Z$ ($n$ factors)
with the lexicographic order. By the unspecified term discrete valuation we always mean
discrete valuation of rank $1$, i.e.\ a discrete valuation in the usual sense.

\ms Let $V$ be a discrete valuation ring of rank $n$ of $K$. Then
$V$ has $n$ distinct prime ideals $\p_0 \supsetneqq \p_1 \supsetneqq \cdots \supsetneqq
\p_n$. For $0 \leq i \leq n$, let $V_i$ be the localization of $V$ at $\p_i$ and let
$k_i=k(\pi_i)$ be the residue field of $V_i$. The $V_i$ are also discrete valuation rings (of
rank $n-i$) of $K$ and
\[
V=V_0 \subset V_1 \subset \cdots V_n=K.
\]
For $0\leq i <j\leq n$ the image of $V_i$ in $k_j$ is a discrete valuation ring of rank
$j-i$. For $0 \leq i < n$, we denote the image of $V_i$ in $k_{i+1}$ by $\bar V_i$. Then
$\bar V_i$ is a usual (i.e.\ rank $1$) discrete valuation ring with quotient field $k_{i+1}$
and residue field $k_i$. By \cite{B-Comm}, Ch.\ VI, \S7, ex.6d, the ring $V$ is henselian if
and only if all $\bar V_i$, $0 \leq i <n$ are henselian (see also \cite{Ri}, ch.\ F, prop.9).

\ms Let $V \subset K$ be a discrete valuation ring of rank $n$
and let $L|K$ be a finite separable extension. Let $w|v$ be a valuation of $L$ extending $v$.
Then also $w$ is a discrete valuation of rank $n$. Assume that the residue characteristic of
$v$ is $p>0$.

\begin{lemma} Let $v$ be a discrete valuation of rank $n$ of a field $K$
and let $L|K$ be a separable extension of $K$. Then $v$ is tamely ramified in $L|K$ if and
only if for every extension $w$ of $v$ to $L$ the ramification index $e_{w|v}$ is prime to
$p$ and  the successive residue extensions of $\bar W_i|\bar V_i$ are separable for $0\leq i
<n$.
\end{lemma}

\demo{Proof:} We may assume that $v$ is strictly henselian and that $L|K$ is finite. In
particular, $v$ has a unique extension $w$ to $L$. Assume that $v$ is tamely ramified in
$L|K$. By proposition \ref{tameprop}, (iii), $L|K$ is an abelian extension of degree prime to
$p$. Since $v$ is defectless in $L|K$, the ramification index is prime to $p$. Furthermore,
the degrees of the successive residue extensions are prime to $p$, and so these extensions
are separable.

It remains to show the other direction. Since the successive residue extensions are
separable, and since discrete valuations of rank $1$ are defectless in separable extensions,
one shows inductively (cf.\ \cite{B-Comm} Ch.VI \S7 ex.5 or \cite{Ri}, G th.5) that $v$ is
defectless in $L|K$. In particular, $[L:K]=e_{w|v}$ is prime to $p$. By proposition
\ref{tameprop} (ii),(iii), we conclude that $L|K$ is contained in a Galois extension of
degree prime to $p$ and is therefore tamely ramified.
\enddemo

 We call a field $k$ an {\bf\boldmath $n$-dimensional henselian field} if it is the
quotient ring of a henselian discrete valuation ring of rank $n$ with finite residue field.
Equivalently, one can give an inductive definition: A $0$-dimensional henselian field is just
a finite field. An $(n+1)$-dimensional henselian field is a field which is henselian under a
discrete (rank $1$) valuation whose residue field is an $n$-dimensional henselian field. So
an $n$-dimensional henselian field $\kappa$ comes along with a sequence
$\kappa_{n-1},\ldots,\kappa_0$ of residue fields, where each $\kappa_i$, $i=n-1,\ldots,0$ is
an $i$-dimensional henselian  field. In particular, $\kappa_0$ is  finite.

We call (by abuse of notation) an extension of $n$-dimensional henselian fields tamely
ramified, if it is tamely ramified with respect to the discrete valuations of rank $n$. The
next result follows easily by induction from the corresponding $1$-dimensional result.

\begin{prop}
Let $\kappa$ be a $n$-dimensional henselian field with last residue field $\kappa_0$ of
characteristic $p >0$ and let $\kappa^t$ be the maximal tamely ramified extension of
$\kappa$. Then there is a natural isomorphism
\[
G(\kappa^t|\kappa)\cong \Big(\; \underset{n\; -times}{\underbrace{\hat \Z^{(p')}(1) \times
\cdots \times \hat \Z^{(p')}(1)}}\; \Big) \rtimes G(\bar{\kappa}_0|\kappa_0),
\]
where $(\hat \Z)^{(p')}$ denotes the prime-to-$p$ part of\/ $\hat \Z$ and $(1)$ means the
first Tate-twist with respect to the cyclotomic character.
\end{prop}

>From now on let $X$ be a reduced, separated, equidimensional scheme of finite type over
$\Spec(\Z)$, in particular (\cite{Ma}, Ch.13, no.34), $X$ is excellent. Let $d=\dim (X)=
\dim_{Krull}(X)$. A {\bf Parshin-chain} $P$ on $X$ is a sequence $(P_0, P_1, \ldots, P_d)$ of
points on $X$ such that
\[
\overline{\{P_0\}} \subset \overline{\{P_1\}} \subset \cdots \subset \overline{\{P_d\}} =X
\]
and $\dim \overline{\{P_i\}}= i$ for $i=0,\ldots, d$.

An inductive localization-henselization procedure, which was proposed by Parshin with
completion instead of henselization and which is described in \cite{K-S}, \S3, provides us
with a functor $(X,P) \longmapsto \kappa_P^h$
\[
\left( \begin{array}{c} \hbox{\rm Reduced schemes $X$ of}\\
             \hbox{\rm finite type over $\Spec(\Z)$}\\
             \hbox{\rm with a Parshin-chain $P$}\\
             \hbox{\rm of length $d=\dim X$}
\end{array}\right)  \lang \left( \begin{array}{c}
                                      \hbox{\rm finite products}\\
                                      \hbox{\rm of $d$-dimensional}\\
                                      \hbox{\rm henselian fields}
                                    \end{array} \right)\, .
\]
Let us briefly recall this construction. First one takes the henselization $\Cal
O^{\,h}_{X,P_0}$ of the local ring of $X$ at $P_0$. Then one considers the finitely many
prime ideals in this ring which lie over $P_1$, passes to the product of the henselizations
of the localizations with respect to this prime ideals, and so on.

\ms
 If in the Parshin-chain $P$ each $P_i$ is a regular point of
$\overline{\{P_{i+1}\}}$, then $\kappa_P^h$ is a $d$-dimensional henselian field rather than
a finite product of such fields.

\ms Assume that $X$ is integral and let $K$ be the function field
of $X$. We say that a discrete valuation ring $V$ of rank $d$ in $K$ {\bf dominates} a
Parshin-chain $P$ if $V_i$ dominates $\Cal O_{X,P_i}$ for $i=1,\ldots, d$. For a proof of the
following proposition see \cite{K-S}, prop.3.3.

\begin{prop} \label{ks3-3.3}
Under the above assumptions,
\[
\kappa_P^h = \prod_V \Quot(V^h),
\]
where $V$ ranges over all discrete valuations of rank $d$ which dominate $P$ and\/
$\Quot(V^h)$ denotes the quotient field of the henselization $V^h$ of $V$.
\end{prop}

\begin{prop} \label{valuecrit}
Let $X$ be a $d$-dimensional, regular, connected scheme of finite type over $\Spec(\Z)$. Let
$U \subset X$ be an open subscheme and $Y=X-U$ the complement. Let $L$ be a finite nilpotent
extension of the function field $K$ of $X$. Then the following are equivalent.

\ms
 \Item{\rm (i)} $U_L \lang U$ is \'{e}tale and tamely ramified along $Y$.\ms
 \Item{\rm (ii)} Every discrete valuation $v$ of rank $d$ in $K$
 which dominates a Parshin chain on $X$  is
 tamely ramified in $L$, and unramified if the dominated chain in contained in
 $U$.\ms
 \Item{\rm (iii)} For every Parshin-chain $P$ on $X$ the extension of
(finite products of) $d$-dimensional henselian fields $\kappa_{P}^hL| \kappa_P^h$ is tamely
ramified, and unramified if $P$ is contained in $U$.
\end{prop}

\demo{Proof:} The implication (i)$\Rightarrow $(ii) is obvious and the equivalence
(ii)$\Leftrightarrow$(iii) follows from proposition \ref{ks3-3.3}. Suppose that (ii) holds.
Let $P_{d-1}$ be a point of codimension $1$ on $U$. We extend $P_{d-1}$ to a Parshin-chain
$(P_0,\ldots, P_{d-1},P_d)$ on $U$. Let $V$ be a discrete valuation ring in $K$ dominating
$P$. Then (notation as above), $V_{d-1}$ dominates and hence is equal to $\Cal
O_{U,P_{d-1}}$. Since $V$ is unramified in $L|K$, so is $V_{d-1}$. By the theorem of
Zariski-Nagata on the purity of the branch locus, we conclude that $U_L \to U$ is \'{e}tale.

It remains to show that $U_L \to U$ is tamely ramified along $Y$. Since $L|K$ is nilpotent,
we can easily reduce to the case that the Galois group is a finite $p$-group for some prime
number $p$. Let $x_1,\ldots,x_n$ be the finitely many codimension~$1$ points which ramify in
$X_L \to X$. Since the set of points in $X$ where $U_L \to U$ is not tame is closed, it
suffices to show that it is tamely ramified at every closed point $y \in Y$.

$U_L \to U$ is tamely ramified at every point of residue characteristic different from $p$
and  at every point which is not contained in $\cup_{i=1}^n \overline{\{x_{i}\}}$.

We show that a closed point $y$ of residue characteristic $p$ and with $y\in
\overline{\{x_{i}\}}$ for some $1 \leq i \leq n$ does not exist. Assume that $y$ is such a
point. Choose a Parshin-chain $P=(P_0,\ldots,P_d)$ on $X$ with $P_0=y$ and $P_{d-1} = x_{i}$
(we find such a chain because $X$ is catenary). Choose a discrete valuation ring $V$ of rank
$d$ in $K$ dominating $P$. Then (for appropriately chosen embeddings to $\bar K$) the
quotient field $K_{x_i}^{\sh}$ of the strict henselization of $\Cal O_{X,{x_i}}$ contains
$\Quot(V^{\sh})$. By assumption $L\Quot(V^{\sh})|\Quot(V^{\sh})$ is a field extension of
degree prime to $p$, hence trivial. But, since $x_i$ is ramified in the  $p$-extension $L|K$,
$LK_{x_i}^{\sh}|K_{x_i}^{\sh}$ is nontrivial. This yields  the required contradiction.
\enddemo

From the proof of the last proposition, we obtain the

\begin{corol} \label{normaltame}
Assume $X$ is connected but only normal. Let $L|K$ be a finite Galois extension of $p$-power
degree. If\/ $U_L\to U$ is \'{e}tale and tamely ramified along $Y$, then $L|K$ is unramified at
every discrete valuation associated to a codimension $1$ point of $X$ which either lies on
$U$ or whose closure in $X$ contains a point of residue characteristic $p$. If $X$ is
regular, the inverse implication is also true.
\end{corol}

\begin{prop} \label{indep2} Suppose that\/ $U$ is regular and let
$X_1$ and $X_2$ be regular schemes which are {\it proper} over $\Spec(\Z)$ and which contain
$U$ as an open subscheme. Let $Y_i=X_i - U$, $i=1,2$, be the closed complements. Let $U' \to
U$ be a finite nilpotent \'{e}tale covering. Then $U' \to U$ is tamely ramified along $Y_1$ if
and only if it is tamely ramified along $Y_2$.
\end{prop}

\demo{Proof:} We may suppose that the degree of $U' \to U$ is a power of some prime number
$p$. Let $X_3$ be a normal scheme which is proper over $\Spec(\Z)$, contains $U$ as an open
subscheme and has proper surjective morphisms $\pi_i:X_3 \to X_1 $, $i=1,2$ making the
diagrams \diagram{ccc}
 U&\injr{}& X_3\\
 \eqd && \mapd{\pi_i}\\
 U& \injr{} &X_i
\enddiagram
commutative. (The existence of such a $X_3$ is well known, for instance, take $X_3$ as the
normalization of the closure of the image of $U$ in $X_1 \times_\Z X_2$.) If $U' \to U$ is
tamely ramified along $Y_1$, then it is tamely ramified along $Y_3=X_3-U$ (see lemma
\ref{tamebasechange}). Let $L$ be the function field of $U'$. Then, by corollary
\ref{normaltame}, $L|K$ is tamely ramified at every discrete valuation on $L$ which is
defined by a codimension $1$ point of $X_3$ which lies on $Y_3$ and whose closure contains a
point of residue characteristic $p$. This set of valuation contains the set of discrete
valuation on $L$ which are defined by a codimension $1$ point of $X_2$ which lies on $Y_2$
and whose closure contains a point of residue characteristic $p$. Applying corollary
\ref{normaltame} again, we conclude that $U' \to U$ is tamely ramified along $Y_2$. \enddemo

Theorem \ref{i1} of the introduction follows immediately from proposition~\ref{indep2}. In
particular, the abelianized tame fundamental group of a regular scheme of finite type over
$\Spec(\Z)$ does not depend on the choice of a compactification. So it is justified to use
the notation $\pi_1^t(X)^{ab}$ for $\pi_1^t(\bar X,\bar X-X)^{ab}$ for this group.

\section{A finiteness result}

The aim of this section is  to prove the following theorem~\ref{p1finite} (=th.\ref{introfin}
from the introduction). In the proof we essentially use a finiteness result of N. Katz and S.
Lang on relative \'{e}tale fundamental groups (\cite{K-L}, th.1).

\begin{theorem} \label{p1finite}
Let ${\Cal O}$ be the ring of integers in a finite extension $k$ of \/ $\mathbb Q$. Let $X$
be a flat $\Cal O$-scheme of finite type whose geometric generic fibre $X \otimes_{\Cal
O}\bar k$ is connected. Assume that $X$ is normal and that the morphism $X \to \Spec (\Cal
O)$ is surjective. Let $U$ be an open subscheme of $X$. Then the abelianized tame fundamental
group $\pi_1^t(X,X-U)^{ab}$ is finite.
\end{theorem}

For the proof of theorem \ref{p1finite} we need the following

\begin{lemma}\label{zpext}
Let $A$ be a strictly henselian discrete valuation ring with perfect (hence algebraically
closed) residue field and with quotient field\/ $k$. Let $k_\infty|k$ be a $\Z_p$-extension.
Let $K|k$ be a regular field extension and let $B\subset K$ be a discrete valuation ring
dominating $A$. Then $B$ is ramified in $Kk_\infty|K$.
\end{lemma}

\demo{Proof:} For each $n \geq 0$, let $k_n|k$ be the unique subextension of degree $p^n$ in
$k_\infty|k$. Let $A_n$ (resp.\ $B_n$) denote the normalization of $A$ (resp.\ $B$) in $k_n$
(resp.\ $Kk_n$). $A_n$ is again a strictly henselian discrete valuation ring. $B_n$ is a
semi-local Dedekind domain.

Suppose that $B$ is unramified in $Kk_\infty$. Fix an $n$. Let $\p$ be the prime ideal of $B$
and let $\p_1,\ldots, \p_g$ be the prime ideals of $B_n$. Since $B_n|B$ is Galois, all $\p_i$
have a common inertia index $f=f_n$, so $N_{B_n|B}(\p_i)=\p^f$ for $i=1, \ldots, g$. Since
$B_n|B$ is \'{e}tale, we have $p^n=[k_n:k]= gf$.

Let $\pi_n$ be a uniformizer of $A_n$. Since $A_n|A$ is totally ramified, we have
\[
v_A(N_{k_n|k}(\pi_n))=1.
\]
Considered as an element of $B_n$, the $\p_i$-valuation of $\pi_n$ is positive and
independent of $i$, $1 \leq i \leq g$. Denoting this positive number by $a$, we have the
following equality of ideals in $B_n$
\[
(\pi_n)= (\p_1 \cdots \p_g)^a.
\]
Hence $N_{B_n|B}((\pi_n))=\p^{afg}$ and therefore
\[
v_B(N_{k_n|k}(\pi_n))= a p^n \geq p^n.
\]
On the other hand, we have
\[
v_B(N_{k_n|k}(\pi_n))= v_B(\pi_0) \cdot v_A(N_{k_n|k}(\pi_n)) = v_B(\pi_0),
\]
where $\pi_0$ is any uniformizer in $A=A_0$. Since $n$ was arbitrary, the inequality
$v_B(\pi_0)\geq p^n$ yields a contradiction.
\enddemo

\demo{Proof of theorem \ref{p1finite}:} Since $X$ is normal, for any open subscheme $V$ of
$U$ the natural homomorphism $\pi_1(V)\to \pi_1(U)$ is surjective. Therefore also the
homomorphism
\[
\pi_1^t(X,X-V)^{ab} \lang \pi_1^t(X,X-U)^{ab}
\]
is surjective and so we may replace $U$ by a suitable open subscheme and assume that $U$ is
smooth over $\bar{S}=\Spec(\Cal O)$. Let $S\subset \bar S$ be the image of $U$. Consider the
commutative diagram
 \diagram{ccccccc}
 0&\lang & \Ker(U/S)& \lang& \pi_1(U)^{ab} & \lang &\pi_1(S)^{ab}\\
 &&\mapd{}&&\surjd{} &&\surjd{}\\
 0&\lang & \Ker^t(U/S)& \lang& \pi_1^t(X,X-U)^{ab} &
 \lang &\pi_1^t(\bar S,\bar S-S)^{ab}
 \enddiagram
where the groups $\Ker(U/S)$ and $\Ker^t(U/S)$ are defined by the exactness of the
corresponding rows, and the two right vertical homomorphisms are surjective. By a theorem of
N. Katz and S. Lang (\cite{K-L}, th.1), the group $\Ker(U/S)$ is finite. By classical
one-dimensional class field theory, the group $\pi_1^t(\bar S,\bar S-S)^{ab}$ is finite (it
is the Galois group of the ray class field of $k$ with modulus $\prod_{\p \notin S}\p$). The
kernel of $\pi_1(S)^{ab}\to \pi_1^t(\bar S,\bar S-S)^{ab}$ is generated by the ramification
groups of the primes of $\bar S$ which are not in $S$. Denoting the product of the residue
characteristics of these primes by $N$, we see that $\pi_1(S)^{ab}$ is the product of a
finite group and a topologically finitely generated pro-$N$ group. Therefore the same is also
true for $\pi_1(U)^{ab}$ and for $\pi_1^t(X,X-U)^{ab}$. By the snake lemma, it will suffice
to show that the cokernel $C$ of the induced map $\Ker(U/S) \to \Ker^t(U/S)$ is a torsion
group.

Let $K$ be the function field of $X$ and let $k_1$ be the maximal abelian extension of $k$
such that the normalization $U_{Kk_1}$ of $U$ in the composite $Kk_1$ is \'{e}tale over $U$ and
tamely ramified along $Y$. By \cite{K-L}, lemma 2, (2), the normalization of $S$ in $k_1$ is
ind-\'{e}tale over $S$. Let $k_2|k$ be the maximal subextension of $k_1|k$ such that the
normalization $S_{k_2}$ of $S$ in $k_2$ is \'{e}tale over $S$ and tamely ramified along $\bar S
-S$. Then $G(k_2|k)=\pi_1^t(\bar S,\bar S-S)^{ab}$ and, by the snake lemma, $C\cong
G(k_1|k_2)$.

In order to show that $C$ is a torsion group, we therefore have to show that $k_1|k_2$ does
not contain a $\Z_p$-extension of $k_2$ for any prime number $p$. Since $k_2|k$ is a finite
extension and $k_1|k$ is abelian, this is equivalent to the assertion that $k_1|k$ contains
no $\Z_p$-extension of $k$ for any prime number $p$. Let $p$ be a prime number and suppose
that $k_\infty|k$ is a $\Z_p$-extension such that the normalization $U_{Kk_\infty}$ is
ind-\'{e}tale over $U$ and ind-tamely ramified along $Y$. A $\Z_p$-extension is unramified
outside $p$ and at least ramified at one prime dividing $p$, see e.g.\ \cite{NSW},
(10.3.20)(ii). Since the normalization of $S$ in $k_1$ is ind-\'{e}tale, we may suppose that
$p|N$. Let $k'$ be the maximal unramified subextension of $k_\infty|k$ and let $\bar S'$ be
the normalization of $\bar S$ in $k'$. Then the base change $X'=X\times_{\bar S}\bar S'\to X$
is \'{e}tale. Hence $X'$ is normal and the pre-image $U'$ of $U$ is smooth and geometrically
connected over $k'$. So, after replacing $k$ by $k'$, we may suppose that $k_\infty|k$ is
totally ramified at a prime $\p |p$, $\p \in \bar S -S$.

Let $\Cal O_\p$ be the local ring of $\bar S$ at $\p$. After a base change to the strict
henselization $A$  of $\Cal O_\p$, we arrive at the following situation (several letters get
a new meaning):

\ms \Item{(1)} a strictly henselian discrete valuation ring $A$
with perfect residue field of characteristic $p$ and with quotient field $k$ of
characteristic $0$ \ms \Item{(2)} a flat connected and normal $A$-scheme of finite type $X$
 which projects surjectively  to $\Spec (A)$ \ms
\Item{(3)} the function field $K$ of $X$ is a regular extension
of $k$
 (i.e.\ $k$ is algebraically closed in $K$) \ms
\Item{(4)} an open subscheme  $U \subset X$ which is smooth over
$K$\ms \Item{(5)} a $\Z_p$-extension $k_\infty|k$ (in which $A$ is totally ramified) \ms
\Item{(6)} The normalization of $U$ in $Kk_\infty$ is ind-\'{e}tale over $U$ and
 ind-tamely ramified along $Y=X-U$.

\ms\noi Let us show that such a situation cannot occur. Let $P$
be the generic point of an irreducible component of the special fibre of $X$ over $A$. Then
$B= \Cal O_{X_P}$ is a discrete valuation ring dominating $A$. Now lemma \ref{zpext} shows
that $B$ ramifies in $Kk_\infty$. Therefore the order of the inertia group of $P$ in
$Kk_\infty$ is divisible by $p=\hbox{\rm char}\,k(P)$. By corollary \ref{tame-inert}, this
contradicts assumption (6).
\enddemo

In the special case $U=X$ we obtain the following corollary.

\begin{corol}
Let ${\Cal O}$ be the ring of integers in a finite extension $k$ of \/ $\mathbb Q$. Let $X$
be a flat $\Cal O$-scheme of finite type whose geometric generic fibre $X \otimes_{\Cal
O}\bar k$ is connected. Assume that $X$ is normal and that the morphism $X \to \Spec (\Cal
O)$ is surjective. Then the abelianized \'{e}tale fundamental group $\pi_1(X)^{ab}$ is finite.
\end{corol}

\vfill \noi{{\sc Alexander Schmidt}, Mathematisches Institut, Universit\"{a}t Heidelberg, Im
Neu\-enheimer Feld 288, 69120 Heidelberg, Deutschland\\ e-mail:
schmidt@mathi.uni-heidelberg.de}

\end{document}